\documentclass[12pt]{article}
\usepackage{amssymb,amsmath,latexsym}
\usepackage{amsbsy}
\usepackage{amsfonts}
\usepackage{amsopn}
\usepackage{amsthm}

\DeclareMathOperator{\p}{\mbox{\rm I\hspace{-0.02in}P}}
\DeclareMathOperator{\e}{\mbox{\rm I\hspace{-0.02in}E}}

\begin{document}

\begin{doublespace}

\newtheorem{thm}{Theorem}[section]
\newtheorem{lemma}[thm]{Lemma}
\newtheorem{defn}{Definition}[section]
\newtheorem{prop}[thm]{Proposition}
\newtheorem{corollary}[thm]{Corollary}
\newtheorem{remark}[thm]{Remark}
\newtheorem{example}[thm]{Example}
\numberwithin{equation}{section}

\newcommand{\D}{{\rm d}}
\def\ee{\varepsilon}
\def\qed{{\hfill $\Box$ \bigskip}}
\def\MM{{\cal M}}
\def\BB{{\cal B}}
\def\LL{{\cal L}}
\def\FF{{\cal F}}
\def\EE{{\cal E}}
\def\QQ{{\cal Q}}

\def\R{{\mathbb R}}
\def\L{{\bf L}}
\def\E{{\mathbb E}}
\def\F{{\bf F}}
\def\P{{\mathbb P}}
\def\N{{\mathbb N}}
\def\eps{\varepsilon}
\def\wh{\widehat}
\def\pf{\noindent{\bf Proof.} }

\title{\Large \bf Convexity and smoothness of scale functions\\
and de Finetti's control problem}
\author{Andreas E. Kyprianou\thanks{Department of Mathematical
Sciences, University of Bath, Claverton Down, Bath, BA2 7AY, U.K.
E-mail: a.kyprianou@bath.ac.uk}, \, V{\'\i }ctor
Rivero\thanks{Department of Mathematical Sciences, University of
Bath, Claverton Down, Bath, BA2 7AY, U.K. and Centro de
Investigaci\'on en Matem\'aticas (CIMAT A.C.), Calle Jalisco s/n,
Col. Valenciana, A. P. 402, C.P. 36000, Guanajuato, Gto. MEXICO.
E-mail: riverovm@gmail.com} \, and Renming Song\thanks{Department of
Mathematics,University of Illinois, 1409 W. Green Street, Urbana, IL
61801, USA . E-mail: rsong@math.uiuc.edu }}
\date{\today}

\maketitle

\begin{abstract}
Under appropriate conditions, we obtain smoothness and convexity properties of $q$-scale functions
for spectrally negative L\'evy processes. Our method appeals directly to very recent developments in the theory of potential analysis of subordinators. As an application of the latter results to scale functions, we are able to continue  the very
recent work of \cite{APP2007} and \cite{Loe}. 
We strengthen their
collective conclusions by showing, amongst other results, that
whenever the L\'evy measure has a  density which is
log convex then for $q>0$ the scale function $W^{(q)}$ is convex on
some half line $(a^*,\infty)$ where $a^*$ is the largest value at
which $W^{(q)\prime}$ attains its global minimum. As a consequence
we deduce that de Finetti's classical actuarial control problem is
solved by a barrier strategy where the barrier is positioned at
height $a^*$.
\end{abstract}



\section{Introduction}
Recently there has been considerable progress in the potential analysis of subordinators, in particular with the identification of a natural class of subordinators known as special subordinators (see for example \cite{SV06, SV07}). At the same time, there has been a growing body of literature
concerning actuarial mathematics which explores the interaction of
classical models of risk and fine properties of L\'evy processes
with a view to gaining new results on both sides (see for example
\cite{APP2007, DK2006, HPSV2004a, HPSV2004b, KKM2004, KK2006,
KP2007, Loe, RZ2007, SV2007}). 

In this paper we shall marry some of these developments together. We will use new potential analytic considerations found in, for example \cite{SV06, SV07}, to understand better smoothness properties of scale functions for spectrally negative L\'evy processes. This builds on other recent developments which closely link the theory of scale functions to potential analysis of subordinators, see \cite{HK2007} and \cite{KR}. In  turn this will 
allow us to  solve de Finetti's classical actuarial control problem for a much larger class of driving spectrally negative L\'evy processes than previously known.  For the
remainder of this introduction we shall elaborate on the latter in
more detail before moving on to our results and their proofs.

Henceforth we assume that $X=(X_t: t\ge 0)$ is a spectrally negative
L\'evy process with L\'evy triplet given by $(\gamma,
\sigma, \Pi)$, where $\gamma \in \R$, $\sigma\ge 0$ and $\Pi$ is a
measure on $(0,\infty)$ satisfying
$$
\int^{\infty}_0(1\wedge x^2)\Pi(dx)<\infty.
$$
The Laplace exponent of $X$ is given by
$$
\psi(\theta)=\log(\E(e^{\theta X_1}))=\gamma
\theta+\frac12\sigma^2\theta^2 -\int^{\infty}_0(1-e^{-\theta
x}-\theta x1_{\{0<x<1\}})\Pi(dx).
$$
(The reader will note that, for convenience, we have arranged the representation of the Laplace exponent in such a way that the support of the L\'evy measure is positive even though the process experiences only negative jumps).
Let $\Phi(0)$ be the largest real root of $\psi$ and recall that $\Phi(0)>0$ if and only if $X$ drifts to $-\infty,$ or equivalently $\psi^{\prime}(0+)<0.$ The restriction $\psi:[\Phi(0),\infty)\to[0,\infty)$ is a bijection whose inverse will be denoted by  $\Phi.$  

Let $\phi$ be the Laplace
exponent of the descending ladder height subordinator
$\widehat{H}=(\widehat{H}_{s}, s\geq 0)$ associated to $X$. Standard
theory  dictates that $\phi$ and $\psi$ are related by the
Wiener-Hopf factorization
$$
\psi(\theta)=(\theta-\Phi(0))\phi(\theta),\qquad \theta\geq 0,
$$
where $\phi$ satisfies
\begin{equation}\label{bf}
\phi(\theta)=\kappa+\D\theta+\int^\infty_{0}
(1-e^{-\theta x})\Upsilon(x)dx,\qquad \theta\geq 0,
\end{equation}
with $\D=\sigma^2/2$, $\kappa\geq 0,$ $\kappa\Phi(0)=0$ and
$\Upsilon:(0,\infty) \to(0,\infty)$ a  function such
that $\int^{\infty}_{0}(1\wedge x) \Upsilon(x)dx<\infty.$  Moreover,
$$
\overline{\Pi}(x):=\int^\infty_{x}\Pi(dx)\text{ and }
\overline{\Upsilon}(x):=\int^\infty_{x}\Upsilon(z)dz =
e^{\Phi(0)x}\int^\infty_{x}e^{-\Phi(0)z}\overline{\Pi}(z)dz, \qquad
x>0
$$
where the last equality is also a well established fact.
The Wiener-Hopf factorization for
$\psi,$ in its Laplace transform form also states that $\psi,$
$\Phi$ and the Laplace exponent of the bivariate descending ladder
processes, say $\widehat{\kappa}: \R^+\times\R^+\mapsto\R$, are
related by the equation
\begin{equation}\label{WH}
\widehat{\kappa}(\alpha,\beta)=c\frac{\alpha-\psi(\beta)}
{\Phi(\alpha)-\beta},\qquad \alpha,\beta\geq 0,
\end{equation}
where $c>0$ is an arbitrary  constant depending on the normalization
of local time at the infiumum. Without loss of generality we can and
will suppose that it is  equal to 1.

A key object in the fluctuation theory of spectrally negative L\'evy
processes and its applications is the {\it scale functions}. For
each $q\geq0$ the so called $q$-scale function of $X$, $W^{(q)}:
\R\to [0, \infty),$ is the unique function such that $W^{(q)}(x)=0$
for $x<0$ and on $[0, \infty)$ is a strictly increasing and
continuous function whose Laplace transform is given by
$$
\int^{\infty}_0e^{-\theta x}W^{(q)}(x)dx=\frac1{\psi(\theta)-q},
\quad \theta>\Phi(q).
$$

In the last 10 years or so the use of scale functions has proved to be of great importance in a wide variety of applied probability models driven by spectrally negative L\'evy processes. We refer to \cite{Kyp}, \cite{HK2007} and \cite{KR} for a recent overview of their presence in the literature. As alluded to above, we are concerned  here in particular with their importance in one of the most classical problems of modern actuarial mathematics: de Finetti's control problem.

Recall
that the classical Cram\'er-Lundberg risk process corresponds to a spectrally negative
L\'evy process  $X$ taking the form of
 a compound Poisson process with arrival rate $\lambda>0$ and
negative jumps, corresponding to claims, having common distribution
function $F$ with finite mean $1/\mu$ as well as a drift $c>0$,
corresponding to a steady income due to premiums. It is usual to
assume the {\it net profit condition} $c - \lambda/\mu >0$ which
says nothing other than $\psi^{\prime}(0+)>0$.

An offshoot of the classical ruin problem for the Cram\'er-Lundberg
process was introduced by de Finetti \cite{DeF}. His intention was
to make the study of ruin under the Cram\'er-Lundberg dynamics more
realistic by introducing the possibility that dividends are paid out
to share holders up to the moment of ruin. Further, the payment of
dividends should be made in such a way as to optimize the expected
net present value of the total dividends paid to the shareholders
from time zero until ruin. Mathematically speaking, de Finetti's
dividend problem amounts to solving a control problem which we state
in the next paragraph but within the framework of the general L\'evy
insurance risk process. The latter process is nothing more than a
general spectrally negative L\'evy process which respects the
analogue of the net profit condition, namely $\psi^{\prime}(0+)>0$ (although
the latter is not necessary in what follows).

Suppose that $X$ is a general spectrally negative L\'evy process (no
assumption is made on its long term behaviour) with probabilities
$\{\mathbb{P}_x : x\in\mathbb{R}\}$ such that under $\mathbb{P}_x$
we have $X_0 = x$ with probability one. (For convenience we shall
write $\mathbb{P}_0= \mathbb{P}$). Let $\xi=\{L^\xi_t: t\geq 0\}$ be
a dividend strategy consisting of a left-continuous non-negative
non-decreasing process adapted to the (completed and right
continuous) filtration $\{\mathcal{F}_t : t\geq 0\}$ of $X$. The
quantity $L^\xi_t$ thus represents the cumulative dividends paid out
up to time $t$ by the insurance company whose risk process is
modelled by $X$. The controlled risk process when taking into
account of the dividend strategy $\xi$ is thus $U^\xi=\{U^\xi_t :
t\geq 0\}$ where $U^\xi_t = X_t - L^\xi_t$. Write $\sigma^\xi
=\inf\{t>0: U^\xi_t < 0 \}$ for the time at which ruin occurs when
the dividend payments are taking into account. A dividend strategy
is called admissible if at any time before ruin a lump sum dividend
payment is smaller than the size of the available reserves; in other
words $L^\xi_{t+} - L^\xi_t \leq\max\{ U^\xi_t,0\}$ for $t\leq\sigma^\xi$.
Denoting the set of all admissible strategies by $\Xi$, the expected
value discounted at rate $q>0$ of the dividend policy $\xi\in\Xi$
with initial capital $x\geq 0$ is given by
$$
v_\xi(x) = \mathbb{E}_x\left(\int_{[0,{\sigma^\xi}]} e^{-qt}
dL_t^\xi\right),
$$
where $\mathbb{E}_x$ denotes expectation with respect to
$\mathbb{P}_x$ and $q>0$ is a fixed rate. {\it De Finetti's dividend
problem} consists of solving the following stochastic control
problem: characterize
\begin{equation}
\label{controlproblem}
v^*(x):= \sup_{\xi\in\Xi} v_\xi(x)
\end{equation}
and, further, if it exists, establish a strategy $\xi^*$ such that
$v^*(x) = v_{\xi^*}(x)$.

This problem was considered by Gerber \cite{G69} who proved that,
for the Cram\'{e}r-Lundberg model with exponentially distributed
jumps, the optimal value function is a result of a barrier strategy.
That is to say, a strategy of the form $L^a_t = a\vee\overline{X}_t
- a$ for some  $a\geq 0$ where $\overline{X}_t :=\sup_{s\leq t}X_s$. In that case the controlled process
$U^{a}_t = X_t - L^a_t$ is a spectrally negative L\'evy process
reflected in the barrier $a$.

This result has been re-considered very recently in \cite{AM2005}
for Cram\'er-Lundberg processes with a general jump distribution. In
the latter paper it was shown that for an appropriate choice of jump
distribution, the above described barrier strategy is not optimal.
In much greater generality, the paper \cite{APP2007} focuses on the
spectrally negative case and finds sufficient conditions for the
optimal strategy to consist of a simple barrier strategy. It is in
the latter paper that we first begin to see the connection with
scale functions as the sufficient conditions given in \cite{APP2007}
are phrased in terms of a variational inequality involving the value
of a barrier strategy which itself can be expressed in terms of the
associated scale function $W^{(q)}$. In a remarkable development
shortly thereafter, Loeffen \cite{Loe} made a decisive statement
connecting the shape of the scale function $W^{(q)}$ to the
existence of an optimal barrier strategy. Loeffen's result begins by
requiring that the scale function $W^{(q)}$ is {\it sufficiently
smooth} meaning that it
 belongs
to $C^1(0,\infty)$ if $X$ is of bounded variation and otherwise
belongs to $C^{2}(0,\infty)$. Loeffen's theorem reads as follows.

\begin{thm}\label{ronnie}
Suppose that $X$ is such that its scale functions are sufficiently
smooth. Let
$$
a^*=\sup\{a\geq 0 : W^{(q)\prime} (a)\leq W^{(q)\prime}(x) \text{
for all }x\geq 0\},
$$
(which is necessarily finite) where we understand $W^{(q)\prime}(0)=
W^{(q)\prime}(0+)$. Then the barrier strategy at $a^*$ is an
optimal strategy if
\begin{equation}
W^{(q)\prime} (a)\leq W^{(q)\prime} (b)\text{ for all }a^*
\leq a\leq b<\infty.
\label{ult-convex}
\end{equation}
\end{thm}

The condition (\ref{ult-convex}) is tantamount to saying that the
scale function $W^{(q)}$ is convex beyond the global minimum of its
first derivative.  An intriguing result in itself, it is however
arguably not a particularly practical condition to verify.
None-the-less \cite{Loe} makes one further striking step by
providing a very natural class of L\'evy risk processes for which
(\ref{ult-convex}) holds. More precisely, it is shown that
(\ref{ult-convex}) holds when the L\'evy measure $\Pi$ is absolutely
continuous with a completely monotone density.

Thanks then to Theorem \ref{ronnie} a clear mandate is set with
regard to finding as broad a class of L\'evy processes as possible
for which the barrier strategy is optimal through smoothness and
convexity properties of the scale functions $W^{(q)}$. Motivated by
this problem this paper serves a twofold purpose. Firstly to
establish results which discuss the issue of smoothness and
convexity of scale functions and secondly, using some of the latter
results, to return to de Finetti's control problem and establish a
larger class of L\'evy processes for which the barrier strategy is
optimal.

The remainder of the paper is structured as follows. In Section
\ref{CSSF} we present an ensemble of results which provide
sufficient conditions for smoothness, concavity and (ultimate)
convexity of scale functions. Key to some of the results in this section 
are  recent potential analytic developments in the theory of subordinators. 
In Section \ref{deF} we give our main
result on de Finetti's control problem: when the L\'evy measure of
the underlying process has a  log convex density,
the solution to de Finetti's control problem is a barrier strategy.
We then make a few remarks about this result and the main issues
involved in the proof of this result. Also in this section we
explain why this is a broader class of L\'evy processes by giving
some explicit examples. In Section \ref{proof} we use the results of
Section \ref{CSSF} to prove our main result on de Finetti's control
problem. We are not able to apply Theorem \ref{ronnie} verbatim for
the present case however. Instead we must revisit its proof in order
to weaken the meaning of `sufficiently smooth' in its statement.
Ultimately this requires the involvement of stochastic calculus
which appeals to both semi-martingale local time and Markov local
time. Some of the proofs of Section \ref{CSSF} are left to an
Appendix.

\section{Convexity and Smoothness of Scale Functions}\label{CSSF}

We will first deal with $0$-scale functions for spectrally negative
L\'evy processes that do not drift to $-\infty,$ that is, processes for which
$\Phi(0)=0.$ Unless otherwise stated throughout this section we will assume that the measure
$\Pi$ has a strictly positive density $\pi(x),$ $x>0,$ with respect
to the Lebesgue measure. In this case,
$$
\Upsilon(x)=\Pi(x,\infty):=\overline{\Pi}(x)=\int^\infty_{x}
\pi(y)dy,\qquad x>0.
$$
Before stating our first result we recall that a subordinator $H$ is said special if there exists another subordinator $H^*,$ the so-called conjugate, such that if $h$ and $h^*$ are their respective Laplace exponents then
$$
\theta=h(\theta)h^*(\theta),\qquad \theta\geq 0.
$$ We refer to \cite{SV07} for a recent account of properties of this subclass of subordinators. We also  mention here in particular that the identification of this class of subordinators has permitted quite significant developments in their potential analysis. Indeed it is the latter developments which play a significant role in the forthcoming analysis of scale functions.

Our first result is on the concavity of the $0$-scale function which, for convenience, we henceforth denote by  $W$ instead of $W^{(0)}$.

\begin{thm}\label{thm:2.1}
Assume that $\Phi(0)=0.$ If the function
$x\mapsto\overline{\Upsilon}(x):=\int^\infty_{x}\Upsilon(z)dz$ is
log convex on $(0, \infty)$, then the scale function $W$ is concave
on $(0, \infty)$.
\end{thm}

\begin{proof} It follows from the log convexity of $\overline{\Upsilon}$ and
Theorem 2.4 of \cite{SV07} that
$\widehat{H}$ is a special subordinator and the renewal function
of $\widehat{H}$ has a decreasing derivative $u$ which is also
called the potential density of $\widehat{H}$. Since $W^{\prime}(x)=u(x)$, we know
that $W$ is concave.
\end{proof}

The next theorem is one of our main results of this section.

\begin{thm}\label{thm:2.2} Assume that $\Phi(0)=0.$ If the function
$x\mapsto \overline{\Pi}(x)= \int^{\infty}_x\pi(s)ds$ is log convex
on $(0, \infty)$, then the function $W^{\prime}$ is convex on $(0, \infty)$.
Furthermore, if $X$ has a Gaussian term or equivalently the drift of
the descending ladder height process is strictly positive then $W\in
C^2(0,\infty).$
\end{thm}

\begin{proof} By our assumption, we know that the function
$$
\overline{\Upsilon}(x)= \int^{\infty}_x\overline{\Pi}(s)ds, \quad x>0
$$ is in $C^1(0,\infty).$ It follows from the first paragraph in
the proof of Theorem 2 in \cite{Gri78} that this function is also
log convex. Therefore it follows from Theorem 2.4 of \cite{SV07}
that $\widehat{H}$ is a special subordinator and the renewal
function of $\widehat{H}$ has a decreasing derivative $u$ which is
also called the potential density of $\widehat{H}$. It follows from
\cite{SV06} that the function $u$ satisfies the following equation
$$
\D u(t)+\int^t_0\overline{\Upsilon}(t-s)u(s)ds=1, \quad t>0,
$$
where $\D\ge 0$ is the drift of $\widehat{H}$. Now when $\D=0$ we
can apply Theorem 3 of \cite{Gri80} to conclude that the function
$u$ is convex.  When $\D>0$ we can apply Theorem 2 of \cite{Gri78}
combined with the first two sentences of Section 4 in \cite{Gri78},
to conclude that $u$ is convex and in $C^1(0,\infty).$ Now the
conclusion follows since $W^{\prime}(x)=u(x)$.
\end{proof}
The two theorems above and the arguments used in their proofs have several consequences, the first of which can be summarized as
follows. If $\overline{\Pi}(s), s>0,$ is a log convex function, then
$\overline{\Upsilon}(s),$ $s>0$, is also log convex and $W^{\prime}$ is a
decreasing and convex function which implies that the subordinator
$\widehat{H}$ with the tail of its  L\'evy measure given by
$\overline{\Upsilon}$ is special, and thus there exists a
subordinator $\widehat{H}^*$, with the tail of its L\'evy measure
denoted by $\overline{\Upsilon}^*(x)$, such that
$$
W^{\prime}(x)=W^{\prime}(\infty)+\overline{\Upsilon}^*(x), \quad
x>0,
$$
and as a consequence $\Upsilon^*$ has a decreasing density in
$(0,\infty).$ Then Theorem 2 and Corollary 1 in \cite{KR} imply in
turn that there exists a spectrally negative L\'evy process that
does not drift to $-\infty$ such that its scale function $W^*$ satisfies$$
W^{*\prime}(x)=\kappa+\overline{\Upsilon}(x), \quad x>0
$$
and therefore $W^{*\prime}$ is log convex.

Another interesting consequence provides a sufficient condition in
terms of the potential density to guarantee that a subordinator has
a L\'evy measure with a decreasing density. This may be useful in
the cases where a subordinator is characterized by its potential
measure, as in the case of subordinators arising in the random covering of the positive reals, see e.g. \cite{FFS}.

\begin{corollary}\label{cor:1}
Let $H$ be a subordinator such that its potential measure has a
density, say $W^{\prime},$ in $(0,\infty)$ such that $W^{\prime}$ is non-increasing
and $-W^{\prime\prime}$ is non-increasing and log convex. Then the L\'evy measure
of $H$ has non-increasing density.
\end{corollary}

\begin{proof}
Since the potential measure of $H$ has a non-increasing density, we
know that $H$ is a special subordinator whose conjugate we will
denote by $H^*.$ Furthermore, the tail of the L\'evy measure of
$H^*$ equals $W^{\prime}(x)-W^{\prime}(\infty), x>0,$ and then its
density is given by $-W^{\prime\prime}.$ We now argue as in Theorem
\ref{thm:2.2} to ensure that the potential measure of $H^*$ admits a
decreasing and convex density in $(0,\infty).$ This finishes the
proof since the tail of the L\'evy measure of $H$ equals the density
in $(0,\infty)$ of the potential measure of $H^*.$
\end{proof}

An interesting question is whether a given function is the scale
function of a spectrally negative L\'evy process. It has been proved
in Corollary 2 in \cite{KR} that a sufficient condition is that such
a function is a Bernstein function. In the next result we provide a
weaker sufficient condition.

\begin{corollary}\label{cor:2} Suppose that $W$ is a function
on $\R$ such that $W(x)=0$ for all $x<0$ and that $W$ is positive
and continuous on $[0, \infty)$. If $W$ is a concave non-decreasing
function on $[0, \infty)$ such that $W^{\prime}$ is non-increasing
on $(0, \infty)$ with $a:=\lim_{x\downarrow 0}xW^{\prime}(x)<\infty$
and $-W^{\prime\prime}$ is non-increasing and log convex on $(0,
\infty)$, then there exists a spectrally negative L\'evy process
such that $W$ is its $0$-scale function.
\end{corollary}

\begin{proof}
We claim that $-W^{\prime\prime}$ is the L\'evy density of some
subordinator, that is,
\begin{equation}\label{c4ld}
-\int^{\infty}_0(1\wedge x)W^{\prime\prime}(x)dx<\infty
\end{equation}
In fact, since $W^{\prime}$ is non-increasing,  we have
$$
-\int^{\infty}_1W^{\prime\prime}(x)dx=W^{\prime}(1)-
W^{\prime}(\infty)<\infty.
$$
On the other hand, since $W^{\prime\prime}$ is non-decreasing, we
have for any $x\in (0, 1)$,
\begin{eqnarray*}
-\int^1_xyW^{\prime\prime}(y)dy&=&-\int^1_xydW^{\prime}(y)=
xW^{\prime}(x)-W^{\prime}(1)+\int^1_xW^{\prime}(y)dy\\
&=&xW^{\prime}(x)-W^{\prime}(1)+W(1)-W(x)\\
&\le&W(1)-W(0)-W^{\prime}(1)+\lim_{x\downarrow 0}xW^{\prime}(x)\\
&\le&W(1)-W(0)-W^{\prime}(1)+a.
\end{eqnarray*}
Thus the claim is valid. 

We may now deduce that $-W^{\prime\prime}$ is the L\'evy density of some special subordinator $H^{*},$
whose conjugate, $H,$ admits $W$ as its potential measure. By Corollary \ref{cor:1} the L\'evy
measure of $H$ has a non-increasing density. Thus we can construct a spectrally negative L\'evy process $X$ whose
descending ladder height process is $H,$ see for example \cite{HK2007}, which also satisfies the assertions
of Corollary 2.4. 
\end{proof}

As we mentioned before the conditions in Corollary~\ref{cor:2} are
weaker than those in Corollary 2 in \cite{KR}, because every
Bernstein function $f$ is a non-decreasing concave function, and
$f^{\prime}$ and $-f^{\prime\prime}$ are completely monotone
functions with $\lim_{x\downarrow 0}xf^{\prime}(x)=0$.

\begin{corollary}\label{cor:subconv}
Let $H$ be a subordinator whose L\'evy density, say $\Upsilon(x), $
$x>0,$ is  log convex (and hence non-increasing) then the restriction of the
potential measure to $(0,\infty)$ has a non-increasing and convex
density. If furthermore, the drift of $H$ is strictly positive then
the density is in $C^1(0,\infty).$
\end{corollary}

\begin{proof}
As above, using the arguments in the first paragraph in the proof of
Theorem 2 in \cite{Gri78} we get that $\overline{\Upsilon}(x)=
\int^\infty_{x}\Upsilon(y)dy,$ $x>0,$ is log convex. Then by Theorem
2.4 in \cite{SV07} we know that $H$ is a special subordinator and
therefore the restriction of its potential measure to $(0,\infty)$
has a non-increasing density. Now, we can simply repeat the
arguments in the proof of Theorem \ref{thm:2.2} to obtain the
convexity and $C^1(0,\infty)$ result.
\end{proof}

The following result is the analogue of Theorems~\ref{thm:2.1} and \ref{thm:2.2} for $q$-scale functions for $q>0$ if $\Phi(0)=0$ and $q\geq 0$ if $\Phi(0)>0.$

\begin{thm}\label{thm:2.3} If the function
$$\overline{\Pi}(x):= \int^{\infty}_x\pi(s)ds, \quad x>0$$
is log convex, then for any $q>0$ if $\Phi(0)=0,$ and $q\geq 0$ if $\Phi(0)>0,$ the function $g_{q}(x):=
e^{-\Phi(q)x}W^{(q)}(x),$ $x>0,$ is concave. If furthermore, the function
$\pi$ is  log convex (and hence non-increasing) then $\overline{\Pi}(x)$ is log convex, the first derivative of
$g_{q}$ is non-increasing and convex and the functions  $W^{(q)}$ and
$W^{(q)\prime}$ are strictly convex in the interval $(a^*,\infty),$
where
$$
a^*=\sup\left\{a\geq 0 : W^{(q)\prime} (a)\leq W^{(q)\prime}(y)
\text{ for all }y\geq 0\right\}<\infty.
$$
Finally, if the latter assumption is satisfied and the Gaussian
coefficient is strictly positive then $W^{(q)}\in C^{2}(0,\infty)$.
\end{thm}

The proof of this theorem relies on the following technical lemmas.
Their proofs will be postponed to the Appendix. Note that in the
first lemma below, the term $q/\Phi(q)$ is to be understood in the
limiting sense, namely $\psi^{\prime}(0+)$, when $q=0$ and $\Phi(0)=0$.

\begin{lemma}\label{teoWH}
For each $q\geq 0,$ the function $\widehat{\kappa}(q,\cdot)$ is a
Bernstein function and its
killing term is given by $$\displaystyle \widehat{\kappa}(q,0)=\frac{q}{\Phi(q)},$$
its drift term is given by  $$\displaystyle \lim_{\theta\to\infty}
\frac{\widehat{\kappa}(q,\theta)}{\theta}=\D$$
and
the tail of its L\'evy measure is given by
 $$\displaystyle \overline{\Upsilon}_{q}(x)
:=e^{\Phi(q)x}\int^\infty_{x}e^{-\Phi(q)y}\overline{\Pi}(y)dy,\qquad x>0.$$
Furthermore, if $\pi$ is non-increasing then for $q\geq 0,$ the L\'evy density associated to $\widehat{\kappa}(q,\cdot)$ is non-increasing.
\end{lemma}

\begin{lemma}\label{logconv}
If $\pi$ is  log convex (and hence non-increasing), then for every $q\geq 0$  the function
$$
\overline{\Pi}(x)-\Phi(q)e^{\Phi(q)x}\int^\infty_{x}e^{-\Phi(q)y}
\overline{\Pi}(y)dy,\qquad x>0,
$$
is log convex.
\end{lemma}

\begin{proof}[Proof of Theorem \ref{thm:2.3}]
We have by assumption that the function $\overline{\Pi}$ is log convex,
which implies that $e^{-\Phi(q)x}\overline{\Pi}(x)$ $x>0$ is also log
convex. Hence it follows from the first paragraph in the proof of
Theorem 2 in \cite{Gri78} that the latter implies that the functions
$$
\int^{\infty}_{x} e^{-\Phi(q)s}\overline{\Pi}(s)ds,\qquad
e^{\Phi(q)x}\int^{\infty}_{x} e^{-\Phi(q)s}\overline{\Pi}(s)ds \quad x>0,
$$
are log convex. It follows that the function
$\overline{\Upsilon}_{q}$ as defined in Lemma~\ref{teoWH} is log
convex and thus by Theorem 2.4 in \cite{SV07}  we have that the
potential density associated to the Bernstein function
$\widehat{\kappa}(q,\cdot)$ has a non-increasing density in
$(0,\infty)$ that we will denote by $u_{q}.$ It follows from Lemmas
1 and 2 in \cite{KR} that the function
$\widehat{\kappa}(q,\Phi(q)+\cdot)$ still is a Bernstein function such that its potential measure  admits the function
 $e^{-\Phi(q)x}u_{q}(x)$ as its
density in $(0,\infty).$ It now follows that the later function is non-increasing
and $\lim_{x\to\infty}e^{-\Phi(q)x}u_{q}(x)=0$.

It is well known that the $q$-scale function $W^{(q)}$ satisfies
$W^{(q)}(x)=e^{\Phi(q)x}W_{\Phi(q)}(x),$ $x>0,$ where $W_{\Phi(q)}$
is the $0$-scale function of the spectrally negative L\'evy process
with Laplace exponent given by
$\psi_{q}(\theta)=\psi(\theta+\Phi(q))-q,$ $\theta\geq 0,$ see e.g.
Lemma 8.4 in \cite{Kyp} for a proof of this fact.
By the Wiener-Hopf
factorization we have that $\psi_{q}$ is given by
$$
\psi_{q}(\theta)=\theta\widehat{\kappa}(q,\Phi(q)+\theta),\qquad
\theta\geq 0.
$$
This implies in turn that
$$
\frac{\theta}{\psi_{q}(\theta)}=\frac{1}{\widehat\kappa(q,\Phi(q) +
\theta)}={\rm d}^*_{q}+\int^\infty_{0}e^{-\theta
x}e^{-\Phi(q)x}u_{q}(x)dx,\qquad \theta\geq 0,
$$
where ${\rm d}^*_{q}=\lim_{\theta\to\infty}1/\widehat{\kappa}
(q,\Phi(q)+\theta)\geq0.$ By the definition of $0$-scale functions
and  integration by parts in the latter equation it follows that
$$
\frac{1}{\psi_{q}(\theta)}=\int^\infty_{0}e^{-\theta x}W_{\Phi(q)}
(x)dx=\int^\infty_{0}e^{-\theta x}\left({\rm d}^*_{q}+\int^x_{0}
e^{-\Phi(q)z}u_{q}(z)dz\right)dx,\qquad \theta\geq 0.
$$
Thus the uniqueness of the Laplace transform implies that
$$
W_{\Phi(q)}(x)={\rm d}^*_{q}+\int^x_{0}e^{-\Phi(q)z}u_{q}(z)dz,
\qquad x\geq0.
$$
Now the first claim immediately follows. The claim about $\overline{\Pi}$ is proved in the proof of Lemma \ref{logconv}. To prove the third claim
we recall that under the assumption that $\pi$ is  log convex (and hence non-increasing)
Lemmas~\ref{teoWH} and \ref{logconv} imply that the L\'evy density
of the Bernstein function $\widehat{\kappa}(q,\cdot)$ is
log convex (and hence non-increasing). Hence the hypotheses of
Corollary~\ref{cor:subconv} are satisfied and therefore $u_{q}$ is a
non-increasing convex function and, whenever the Gaussian
coefficient, equivalently the linear term in $\widehat{\kappa},$ is
strictly positive we have $u_{q}\in C^{1}(0,\infty).$ By elementary
arguments it follows that $e^{-\Phi(q)x}u_{q}(x),$ $x>0,$ satisfies
the same properties. Thus $g_{q}$ is a concave function whose first
derivative is convex and continuous in $(0,\infty)$.

To prove the claim about the convexity of $W^{(q)}$ and $W^{(q)\prime},$
we observe that as $W^{(q)\prime}$ is given by
\begin{equation}\label{eq:qscaleder}
W^{(q)\prime}(x)=\Phi(q)W^{(q)}(x)+u_{q}(x),\qquad x>0;
\end{equation}
and since $u_{q}$ is convex, we will automatically get that
$W^{(q)\prime}$ is ultimately convex once we have proved that $W^{(q)}$
is ultimately convex. Indeed, we have that
\begin{equation}\label{eq:w'}
\begin{split}
W^{(q)\prime\prime}(x)&=\left(\Phi(q)\right)^2W^{(q)}(x)+\Phi(q)u_{q}(x)+u^{\prime}_{q}(x),
\qquad  \text{a.e.}\ x>0.
\end{split}
\end{equation}
Then  as $u^{\prime}_{q}$ increases and $W^{(q)}$ grows exponentially fast
it follows that ultimately $W^{(q)\prime\prime}>0.$ Hence $W^{(q)}$
and $W^{(q)\prime}$ are ultimately strictly convex. Furthermore,
because $W^{(q)\prime}$ tends to infinity as $x$ tends to $\infty,$
it follows that $a^*<\infty.$ Now, let $\alpha_{1}<\alpha_{2}$ be
points at which $W^{(q)\prime}$ reaches a local minimum. Because of
the convexity of $u_{q}$ we know that the right and left derivatives
of $u_{q}$ exist everywhere and they satisfy that
$u^{\prime -}_{q}(\alpha_{1})\leq u^{\prime +}_{q}(\alpha_{1})\leq
u^{\prime -}_{q}(\alpha_{2})\leq u^{\prime +}_{q}(\alpha_{2}).$ As a consequence
the right and left derivatives of $W^{(q)\prime}$ exist everywhere
and satisfy
\begin{equation*}
W^{(q)\prime\prime -}(\alpha_{i})=\Phi(q)W^{(q)\prime}(\alpha_{i})+u^{\prime -}_{q}
(\alpha_{i})\leq 0,\qquad W^{(q)\prime\prime +}(\alpha_{i})=\Phi(q)W^{(q)\prime}
(\alpha_{i})+u^{\prime +}_{q}(\alpha_{i})\geq 0,
\end{equation*}
for $i=1,2.$ These facts together imply that
$$0\leq \Phi(q)\left(W^{(q)\prime}(\alpha_{1})-W^{(q)\prime}(\alpha_{2})\right)
+u^{\prime +}_{q}(\alpha_{1})-u^{\prime -}_{q}(\alpha_{2}),
$$
and hence $W^{(q)\prime}(\alpha_{1})-W^{(q)\prime}(\alpha_{2})\geq
0.$ This implies that the last place where $W^{(q)\prime}$ reaches a
local minimum is also the last place where it hits its global
minimum. Moreover, for $x>0$ we have that $W^{(q)\prime}(a^*)\leq
W^{(q)\prime}(x)$. It thus follows that the following inequalities
\begin{equation*}
\begin{split}
0&\leq \Phi(q)W^{(q)\prime}(a^*)+u^{\prime +}_{q}(a^*)\leq\Phi(q)W^{(q)\prime}(x)+
u^{\prime -}_{q}(x)\leq \Phi(q)W^{(q)\prime}(x)+u^{\prime +}_{q}(x),
\end{split}
\end{equation*}
hold for $x>a^*.$ Actually, the second inequality is a strict one.
Indeed, if there would exist $x^*>a^*$ such that
$\Phi(q)W^{(q)\prime}(a^*) + u^{\prime
+}_{q}(a^*)=\Phi(q)W^{(q)\prime}(x^*)+ u^{\prime -}_{q}(x^*),$ then
since $u^{\prime -}_{q}(x^*)-u^{\prime +}_{q}(a^*)\geq 0,$ we would
have that $W^{(q)\prime}(a^*)\geq W^{(q)\prime}(x^*),$ which would
be a contradiction to the fact that $a^*$ is the largest value where
$W^{(q)\prime}$ attains its global minimum. It follows that
$W^{(q)\prime}$ is strictly increasing for $x>a^*.$ That is
$W^{(q)}$ is strictly convex in $(a^*,\infty)$ and, from equation
(\ref{eq:qscaleder}), we deduce that $W^{(q)\prime}$ is also
strictly convex for $x>a^*$.
Finally, the equation (\ref{eq:w'}) proves also that when
furthermore the Gaussian coefficient is strictly positive then
$W^{(q)}\in C^{2}(0,\infty)$ as in this case we already proved
that $u_{q}\in C^1(0,\infty).$
\end{proof}

We now leave behind the assumption that $\Pi$ has a strictly positive density and allow $\Pi$ to be any L\'evy measure. The following result provides necessary and sufficient conditions
for a scale function associated to a spectrally negative L\'evy
process of bounded variation to be in the class $C^{1}(0,\infty).$

\begin{thm}\label{thm:2.4}
Assume that $X$ is a spectrally negative L\'evy process of bounded
variation. The following conditions are equivalent
\begin{itemize}
\item[(i)] $W\in C^{1}(0,\infty);$
\item[(ii)] $W^{(q)}\in C^{1}(0,\infty)$ for all $q\geq 0;$
\item[(iii)] $\overline{\Pi}\in C^{0}(0,\infty).$
\end{itemize}
\end{thm}

\begin{proof}
First we prove that (i) and (iii) are equivalent. As $X$ is assumed
to be of bounded variation it follows that $X_t=\delta t-S_{t},$
where $\delta>0$ and $S$ is a subordinator. Let $\overline{n}$ be
the excursion measure of $X$ reflected at its supremum. For
background on the excursion theory of L\'evy processes reflected at
their supremum,  see e.g. \cite{Be1996}, \cite{Don2007}, \cite{Kyp}.
A standard result, see e.g. Lemma 8.2 in \cite{Kyp}, says that $W\in
C^1(0,\infty)$ if and only if the law of the height of the excursion
process has no atoms, that is, the measure
$$
\overline{n}\left(\sup_{0\leq s\leq \zeta} \epsilon(s)\in
dx\right), \qquad x>0
$$
has no atoms; where $\epsilon$ and $\zeta$ denotes the generic
excursion process and its lifetime respectively. The proof of the
equivalence of (i) and (iii) will be obtained as a consequence of
the fact above  and the identity
\begin{equation}\label{eq:supexc}
\overline{n}\left(\sup_{0\leq s\leq \zeta}\epsilon(s)>z\right)=
\frac{1}{\delta}\Pi(z,\infty)+\frac{1}{\delta}\int^z_{0}\Pi(dx)
\left(1-\frac{W(z-x)}{W(z)}\right), \qquad z>0.
\end{equation}
If we take this identity for granted, then we have that
\begin{equation}
\begin{split}
\overline{n}\left(\sup_{0\leq s\leq z}\epsilon(s)=z\right)=
\frac{1}{\delta}\Pi\{z\}+\frac{1}{\delta}\Pi\{z\}
\left(1-\frac{W(0)}{W(z)}\right),\qquad z>0.
\end{split}
\end{equation}
It follows then that $\overline{n}(\sup_{0\leq s\leq
\zeta}\epsilon(s)\in dx)$ has atoms if and only if $\Pi$ does. This
proves the equivalence of (i) and (iii).

We will now prove the identity (\ref{eq:supexc}). It is known from
\cite{Rog} and Proposition 5 in \cite{Win} that when $X$ is of
bounded variation the excursion measure of $X$ reflected at its
supremum can be described by the formula
\begin{equation*}
\overline{n}\left(F\left(\epsilon(s), 0\leq s\leq \zeta\right)
\right)=\frac{1}{\delta}\int^\infty_{0}\Pi(dx)\widehat{\e}_{x}
\left(F(X_{s}, 0\leq s\leq \tau^-_{0})\right),
\end{equation*}
where  $F$ is any nonnegative measurable functional on the space of
cadlag paths,  $\widehat{\e}_{x}$ denotes the law of the dual L\'evy
process $\widehat{X}=-X$ and $\tau^-_{x}=\inf\{s>0: X_{s}< x\},$
$x\in\R.$ We will also denote by $\tau^+_{z}=\inf\{s>0: X_{s}>z\},$
$z\in\R.$ Hence, it follows that
\begin{equation*}
\begin{split}
\overline{n}\left(\sup_{0\leq s\leq \zeta}\epsilon(s)>z\right)&=
\frac{1}{\delta}\int^\infty_{0}\Pi(dx)\widehat{\p}_{x}\left(
\sup_{0\leq s\leq \tau^-_{0}}X_{s}>z\right)\\
&=\frac{1}{\delta}\Pi(z,\infty)+\frac{1}{\delta}\int^z_{0}
\Pi(dx)\widehat{\p}_{x}\left(\tau^{+}_{z}<\tau^-_{0}\right)\\
&=\frac{1}{\delta}\Pi(z,\infty)+\frac{1}{\delta}\int^z_{0}
\Pi(dx)\widehat{\p}\left(\tau^{+}_{z-x}<\tau^-_{-x}\right)\\
&=\frac{1}{\delta}\Pi(z,\infty)+\frac{1}{\delta}\int^z_{0}
\Pi(dx)\p\left(\tau^{-}_{x-z}<\tau^+_{x}\right)\\
&=\frac{1}{\delta}\Pi(z,\infty)+\frac{1}{\delta}\int^z_{0}
\Pi(dx)\left(1-\frac{W(z-x)}{W(z)}\right),
\end{split}
\end{equation*}
where in the last equality we have used Takac's solution to the two
sided exit problem for spectrally negative L\'evy processes, see
e.g. \cite{Be1996} Theorem VII.8.

To complete the proof we show that (ii) and (iii) are equivalent.
There is very little work to do as soon as one recalls the fact
stated earlier that for $q>0$
$$
W^{(q)}(x) = e^{\Phi(q)x}W_{\Phi(q)}(x),
$$
where $W_{\Phi(q)}$ is the 0-scale function of the spectrally
negative L\'evy process whose Laplace exponent is given by
$\psi(\theta+\Phi(q))-q$. The latter process has a L\'evy measure
$\Pi_{\Phi(q)}$ given by $\Pi_{\Phi(q)}(dx)=e^{-\Phi(q)x} \Pi(dx)$
on $(0,\infty)$. It follows that $\Pi$ has no atoms if and only if
$\Pi_{\Phi(q)}$ has no atoms. Consequently the arguments above
leading to the equivalence of (i) and (iii) show that $W_{\Phi(q)}$
belongs to $C^1(0,\infty)$ or equivalently (ii) holds, if and only
if (iii) holds.
\end{proof}

\noindent Note that another proof of the last result may be implicitly extracted from Lemma 1 (iii) in \cite{Doney}. In essence that would again necessitate the observation that the entrance law of excursions begin with a jump whose intensity is given by $\delta^{-1}\Pi$.

\section{De Finetti's control problem}\label{deF}

In this section we shall discuss an important consequence of Theorem
\ref{thm:2.3} pertaining  to de Finetti's control problem.
 In
particular we shall prove the following result.

\begin{thm}\label{barrier}
Suppose that $X$ has a L\'evy density $\pi$ that is  log convex
then the barrier strategy at $a^*$ is optimal for (\ref{controlproblem}).
\end{thm}

Before proceeding to the proof of Theorem \ref{barrier}, let us
first make some remarks.

\begin{enumerate}

\item In principle the proof of Theorem \ref{barrier} follows
directly from Theorem \ref{thm:2.3} and Theorem \ref{ronnie} if one
can verify that $W^{(q)}$ is sufficiently smooth. This is possible
in most cases, but not all. The outstanding case is the focus of the
proof of Theorem \ref{barrier} and we identify it below by excluding
the cases for which sufficient smoothness can be established.

If $a^*=0$ then necessarily, either $\sigma>0$, or $\sigma=0$ and
$\Pi(0,\infty)<\infty$ simultaneously. Other types of spectrally
negative L\'evy processes are not possible when $a^*=0$ since then
necessarily $W^{(q)}(0+)=\infty$. In the  case $\sigma>0$ we see
that $W^{(q)}\in C^2(0,\infty)$ by Theorem \ref{thm:2.3} and when
$\sigma=0$ and $\Pi(0,\infty)<\infty$ simultaneously we see from
(\ref{eq:qscaleder}) $W^{(q)}\in C^1(0,\infty)$. Thus when $a^*=0$
we have that $W^{(q)}$ is sufficiently smooth.

Suppose now that $a^*>0$. If $X$ is of bounded variation or
$\sigma>0$ then, similar to the previous paragraph we may again
deduce from Theorem \ref{thm:2.3} and (\ref{eq:qscaleder}) that
$W^{(q)}$ is sufficiently smooth.

The outstanding case is thus given by $a^*>0$, $X$ is of unbounded
variation and $\sigma=0$.

\item Recall that Theorem 3 of \cite{Loe} states that if $X$ has
L\'evy density $\pi$ which is completely monotone then $W^{(q)\prime}$
is convex on $(0,\infty)$ and hence the barrier strategy at $a^*$ is
an optimal strategy. Theorem \ref{barrier} is an improvement on this
result on account of the fact that any completely monotone function
density is  log convex. Below are two some
examples of L\'evy densities which meet the criteria of Theorem
\ref{barrier} but not Theorem 3 of \cite{Loe}.

Suppose that $f$ and $g$ both map $(0,\infty)$ to $[0,\infty)$ and
that they are both non-increasing and log convex. Suppose moreover
that for some (and hence every) $\varepsilon>0$, $\int_0^\varepsilon
x^2 f(x)dx<\infty$ and $\int_\varepsilon^\infty g(x)dx<\infty$.
Further, for some fixed $\alpha>0$ we have $f(\alpha)=g(\alpha)$ and
$$
\frac{f^{\prime -}(\alpha)}{f(\alpha)} \leq \frac{g^{\prime +}
(\alpha)}{g(\alpha)}.
$$
Then
$$
\pi(x) :=\left\{
\begin{array}{ll}
f(x) & x\in(0,\alpha)\\
g(x) & x\in[\alpha,\infty)
\end{array}
\right.
$$
is an example of a decreasing, log convex function which is not
completely monotone in general. Specific cases in which $\pi$ is not
completely monotone may be taken to be
\begin{description}
\item[(i)] $f(x) = x^{-(1+\lambda_1)}$,  $g(x) = x^{-(1+\lambda_2)}$,
$\alpha=1$ where $0<\lambda_2<\lambda_1<2$,
\item[(ii)] $f(x) = e^{2-x}$, $g(x)= e^{1-\lambda x}$, $\alpha
=1/(1-\lambda) $ where $0<\lambda <1$.
\end{description}

\item It is also worth noting that if the L\'evy density $\pi$ meets
the conditions of Theorem \ref{barrier} but is not completely monotone
as in Theorem 3 of \cite{Loe}, then the behaviour of the scale function
$W^{(q)}$ on $(0,a^*)$ is not necessarily concave as is the case in the
aforementioned theorem.

\item The proof of Theorem \ref{barrier}, given in the next section,
is lengthy requiring some auxiliary results first. Scanning the proof
it is not immediately clear where the need for convexity on $(a^*,\infty)$
is needed. The precise point at which this property is required is
embedded in the proof of Lemma \ref{l3} below and we have indicated
as such in the proof.
\end{enumerate}

\section{Proof of Theorem \ref{barrier}}\label{proof}

Following the first remark in the previous section, we shall assume
throughout this section that $a^*>0$, $X$ is of unbounded variation
and $\sigma=0$. Moreover we shall assume that the conditions of Theorem \ref{barrier} are in force.

We define an operator $(\Gamma, {\cal D}(\Gamma))$
as follows. ${\cal D}(\Gamma)$ stands for the family of functions
$f\in C^1(0, \infty)$ such that the integral
$$
\int_{(0,\infty)}[f(x-y) -f(x) + yf^{\prime}(x)\mathbf{1}_{\{y\leq
1\}}]\Pi(dy)
$$
is absolutely convergent for all $x>0$. For any $f\in{\cal
D}(\Gamma)$, we define
$$
\Gamma f(x) = \gamma f^{\prime}(x) + \int_{(0,\infty)}[f(x-y) -f(x) +
yf^{\prime}(x)\mathbf{1}_{\{y\leq 1\}}]\Pi(dy).
$$

Recall that for any $a>0$, the expected value discounted at rate
$q>0$ of the barrier strategy at level $a$ is given by
\begin{equation*}
 v_a(x):= \mathbb{E}_x\left(\int_{[0,\sigma^a]} e^{-qt}
dL^a_t\right) = \left\{
\begin{array}{ll}
W^{(q)}(x)/W^{(q)\prime}(a),
 & -\infty< x\leq a,\\
x-a +W^{(q)}(a)/W^{(q)\prime}(a), &\infty> x>a.
\end{array}
\right.
\end{equation*}
where $\sigma^a=\inf\{t>0: U^a_t<0\}$. The second equality is taken
from \cite{APP2007}.

\begin{lemma}\label{l1} For any $a>0$, $v_a\in{\cal
D}(\Gamma)$. Furthermore, the function $x\mapsto \Gamma v_a(x)$ is
continuous in $(0, a)$.
\end{lemma}

\begin{proof} We have proved in Section 2 that $W^{(q)}$ is in
$C^1(0, \infty)$, hence we know that $v_a$ is in $C^1(0, \infty)$.
To show that $v_a\in{\cal D}(\Gamma)$, we only need to show that the
integral in the definition of $\Gamma v_a$ is absolutely convergent
for all $x>0$. It is easy to check that this is true for $x>a$, so
we are going to concentrate on $x\in (0, a)$. Note that it suffices
to consider $W^{(q)}$ instead of $v_a$. For each $x\in (0,a)$ we may
write the integral in the definition of $\Gamma W^{(q)}$ as
\begin{eqnarray}
&& \int_{(\varepsilon,
\infty)}(W^{(q)}(x-y)-W^{(q)}(x)+yW^{(q)\prime}
(x)\mathbf{1}_{\{y\leq 1\}})\Pi(dy)\notag\\
&&+ \int_{(0,\varepsilon)}(W^{(q)}(x-y)-W^{(q)}(x)+y
W^{(q)\prime}(x))\Pi(dy)\label{abscgt}
\end{eqnarray}
where the value of  $\varepsilon=\varepsilon(x)\in(0,1)$ is chosen
for each $x$ such that $x-2\varepsilon>0$. The absolute convergence
of the first integral as well as its continuity  in $x$ follows in a
straightforward way as a consequence of the continuity and
boundedness of $W^{(q)\prime}$ on bounded intervals and dominated
convergence in the case of continuity. With regard to the second
integral, recall that $W^{(q)\prime}(x)= \Phi(q)W^{(q)}(x) +
u_q(x)$. Using the mean value theorem and the fact that $u_q$ is
convex and decreasing, we get that for  all $y\in(0,\varepsilon)$
\begin{eqnarray*}
&&|W^{(q)}(x-y)-W^{(q)}(x)+yW^{(q)\prime}(x)|\\
&&=y|W^{(q)\prime}(x) -W^{(q)\prime}(x-\xi(y))
|\quad\text{ where }\xi(y)\in(0,y)\\
&&\leq\Phi(q)y|W^{(q)}(x) -W^{(q)}(x-\xi(y)) |+
y|u_q(x) -u_q(x-\xi(y)) |\\
&&\leq \Phi(q)y^2\sup_{z\in[-\varepsilon,\varepsilon]}
W^{(q)\prime}(x+z) +y \int_{x-\xi(y)}^x |u_q^{\prime}(y)|dy\\
&&\leq  \Phi(q)y^2\sup_{z\in[-\varepsilon,\varepsilon]}
W^{(q)\prime}(x+z) +y^2 |u_q^{\prime}(x-\varepsilon)|\\
&&\leq y^2\sup_{z\in[-\varepsilon,\varepsilon]}(\Phi(q)
W^{(q)\prime}(x+z) + |u^{\prime}_q(x+z)|).
\end{eqnarray*}
This estimate shows both that the second integral is uniformly
integrable in (\ref{abscgt}) and continuous in $x$ by dominated
convergence.
\end{proof}

\begin{lemma}\label{l2} For any $a>0$ we have
$$
(\Gamma-q)v_a(x)=0, \quad x\in (0, a).
$$
\end{lemma}

\begin{proof}
It is well known that $e^{-q(t\wedge
\tau^+_{a}\wedge\tau^-_0)}W^{(q)}(X_{t\wedge
\tau^+_{a}\wedge\tau^-_0})$ is a $\mathbb{P}_x$-martingale for each
$x\in(0,a)$ (cf. \cite{AKP2004}), thus $e^{-q(t\wedge
\tau^+_{a}\wedge\tau^-_0)}v_a(X_{t\wedge
\tau^+_{a}\wedge\tau^-_0})$ is a $\mathbb{P}_x$-martingale for each
$x\in(0,a)$.  Appealing to the Meyer-It\^o formula (cf. Theorem 70
of \cite{protter}) 
we have on
$\{t<\tau^+_{a}\wedge\tau^-_{1/n}\}$
\begin{equation}
v_a(X_t) - v_a(x)= m_t+\int_0^t \Gamma v_a(X_{s})ds +\int_\mathbb{R} v_a''(y)\ell(y,t) dy
\qquad \mathbb{P}_x\mbox{-a.s.}
\label{smgw}
\end{equation}
where $\ell(y,\cdot)$ is
the {\it semi-martingale local time} at $y$ of $X$ and, with $X^{(1)}$ as  the
martingale part of $X$ consisting of compensated jumps of size less
than or equal to unity,
\begin{eqnarray*}
m_t &=& \sum_{s\leq t}[\Delta
v_a(X_s)  -\Delta X_s v_a^{\prime}(X_{s-})
\mathbf{1}_{\{|\Delta X_s| \leq 1\}} ] \\
&&- \int_0^t \int_{(0,\infty)}[v_a(X_{s-}-y) -
v_a(X_{s-}) +
yv_a^{\prime}(X_{s-})\mathbf{1}_{\{y\leq
1\}}]\Pi(dy)ds\\
&& + \int_0^t v_a^{\prime}(X_{s-})dX^{(1)}_s
\end{eqnarray*}
is a local martingale which is also a true martingale on account of the fact that  $W^{(q)\prime}$ is
bounded on $[1/n, a]$.  Note that we have used   that the integral part of  $\Gamma v_a(y)$ is
absolutely convergent for each $y\in(0,a)$ in order to meaningfully write down the compensation in the expression for the martingale $m_t$.
The occupation formula
for the semimartingale local time  of $X$ says that
$$
 \int_{\mathbb{R}}\ell(y,t) g(y)dy = \sigma^2\int_0^t
 g(X_{s})ds \qquad \mathbb{P}_x\mbox{-a.s.}
$$
where $g$ is a bounded Borel measurable function. This implies that for Lebesgue almost every $y$, $\ell(y,\cdot)$
is identically zero almost surely. Taking this into account the last
integral in (\ref{smgw}) is almost surely zero.
Using the semi-martingale decomposition in (\ref{smgw}), one may now use stochastic integration by parts for semi-martingales to deduce that on $\{t<\tau^+_a\wedge \tau^-_{1/n}\}$
\[
e^{-qt}v_a(X_t) - v_a(x)= \lambda_t+\int_0^t e^{-qs} (\Gamma -q)v_a(X_{s})ds \qquad \mathbb{P}_x\mbox{-a.s.}
\]
where $\lambda_t = \int_0^t e^{-qs}dm_s$ is a martingale.

Next, use uniform boundedness of $(\Gamma-q)v_a(x)$ on $(\alpha, \beta)\subset[0,a]$ and the martingale property of $e^{-q(t\wedge\tau^+_\beta\wedge\tau^-_\alpha)}v_a(X_{t\wedge\tau^+_\beta\wedge\tau^-_\alpha})$ to get
\begin{equation}
0=\mathbb{E}_x\left[\int_0^{\tau^+_{\beta}\wedge \tau^-_{\alpha}} e^{-qs}(\Gamma - q)v_a(X_{s})ds \right]
=
\int_{(\alpha, \beta)}(\Gamma - q)v_a(y) u^{(q)}(x,y)dy
\label{violate}
\end{equation}
where $u^{(q)}(x,y)dy = \int_0^\infty e^{-qs}\mathbb{P}_x(X_s\in dy; t<\tau^+_\beta\wedge \tau^-_\alpha)$ is the strictly positive resolvent density of the process $X$ killed on exiting $(\alpha, \beta)$ (see \cite{Be1997} for more details). As $(\alpha, \beta)$ is arbitrary and $(\Gamma - q)v_a(x)$ is continuous, it follows that the latter is identically zero on $(0,a)$. This is due to a classical argument by contradiction. If the claim is false then there by continuity of $(\Gamma - q)v_a$ and strict positivity of $u^{(q)}$, there exists an interval $(\alpha',\beta')\subset[0,a]$ such that (without loss of generality)  $(\Gamma-q)v_a(x)>0$ on $(\alpha',\beta')$. Then the equality (\ref{violate}) would be violated.

\end{proof}

For convenience, we use $v$ to denote the function $v_{a^*}$, $U$ to
denote $U^{a^*}$ and $L$ to denote $L^{a^*}$. Then we have the
following result.

\begin{lemma}\label{l3} For any $x>0$ we have
$(\Gamma-q)v(x)\le0.$
\end{lemma}

\begin{proof} There is nothing to prove when $x\in(0,a^*)$ because
of Lemma \ref{l2} applied to the case $a=a^*$. Thanks to the
continuity given by Lemma \ref{l1}, this maybe extended to
$(0,a^*]$. Finally the inequality can be proved to hold on
$(a^*,\infty)$ by following verbatim the arguments in the proof of
Theorem 2 in \cite{Loe}, although it is not necessary to replicate
the behaviour of second derivatives in that proof, since we have
$\sigma=0$.

It is important to note that the use of the convexity of $W^{(q)}$
on $(a^*,\infty)$ appears in Theorem 2 of \cite{Loe} and therefore
in this paper the use of convexity is hidden in the latter part of
the proof. \end{proof}

Now we are in a position to prove Theorem \ref{barrier}.

\begin{proof}[Proof of Theorem \ref{barrier}]
Recall that we are assuming that $a^*>0$, $X$ is of unbounded
variation and $\sigma=0$. Also recall we use $v$ to denote the function $v_{a^*}$, $U$ to
denote $U^{a^*}$ and $L$ to denote $L^{a^*}$. The idea of the proof are similar to that
of \cite{APP2007} and \cite{Loe}, however it is necessary to revisit
the main line of reasoning and provide  more sensitive arguments
that accommodate for the fact that in the present case $W^{(q)}$ is
not {\it sufficiently smooth}, it is twice continuously
differentiable almost everywhere but is not in $C^{2}(0,\infty)$.

Now let $\xi$ be an admissible strategy and let
$\overline{L}^\xi$ be the cadlag modification of the process
$L^\xi$. Note that it is still adapted  as the usual conditions have been assumed on the natural filtration and hence $\overline{U}^\xi = X-
\overline{L}^\xi$  is a semi-martingale.
Since the integral in the definition of $\Gamma v$ is absolutely
convergent for every $x>0$ and second derivative of $v$ is well
defined Lebesgue almost everywhere, we can apply the Meyer-It\^o
formula (cf. Theorem 70 of \cite{protter}) to the process
$v(\overline{U}^\xi )$ 
to get,  after some straightforward manipulation, that on $\{t<\sigma^\xi\}$,
\begin{eqnarray}
v(\overline{U}^\xi_t)-v(\overline{U}^\xi_0) &=& M^\xi_t + \int_0^t
\Gamma v(\overline{U}^\xi_{s-})ds\notag\\
&& +
\sum_{s\leq t}\mathbf{1}_{\{\Delta \overline{L}^{\xi}_s>0\}}
\left\{ v(\overline{U}^\xi_{s-} + \Delta X_s - \Delta \overline{L}^\xi_s) - v(\overline{U}^\xi_{s-} +\Delta X_s) \right\}\notag\\
&&- \int_{(0,t]} v^{\prime}(\overline{U}^\xi_{s-})d\overline{L}^{\xi, c}_s +
\frac{1}{2}\int_\mathbb{R} v^{\prime\prime}(x)\ell^{\xi}(x,t) dx\label{MI}
\end{eqnarray}
such that
\begin{eqnarray*}
M^\xi_t &=& \sum_{s\leq t} \mathbf{1}_{\{|\Delta X_s|>0\}}[
v(\overline{U}^\xi_{s-} + \Delta X_s) - v(\overline{U}_{s-}^\xi)  -\Delta X_s v^{\prime}(\overline{U}^\xi_{s-})
\mathbf{1}_{\{|\Delta X_s| \leq 1\}} ] \\
&&- \int_0^t \int_{(0,\infty)}[v(\overline{U}^\xi_{s-}-y) -
v(\overline{U}^\xi_{s-}) +
yv^{\prime}(\overline{U}^\xi_{s-})\mathbf{1}_{\{y\leq
1\}}]\Pi(dy)ds\\
&& + \int_0^t v^{\prime}(\overline{U}^{\xi}_{s-})dX^{(1)}_s
\end{eqnarray*}
is a local martingale  with $M^{\xi}_0=0$ where $X^{(1)}$ the
martingale part of $X$ consisting of compensated jumps of size less
than or equal to unity. Moreover $\overline{L}^{\xi, c}$ is the
continuous part of $\overline{L}^{\xi}$ and $\ell^\xi(x,\cdot)$ is
the {\it semi-martingale local time} at $x$ of $\overline{U}^\xi$.
We have used in particular the absolute convergence of the integral
part of $\Gamma v$ in order to make sense of $M^\xi_t$ as a
compensated stochastic integral. Similarly to before, the occupation formula
for the semimartingale local time  of $\overline{U}^\xi$ reads
$$
 \int_{\mathbb{R}}\ell^\xi(x,t) g(x)dx = \sigma^2\int_0^t
 g(\overline{U}^\xi_{s})ds,
$$
where $g$ is a bounded Borel measurable function. Also similarly to before, since $\sigma=0$
this implies that for Lebesgue almost every $x$, $\ell^\xi(x,\cdot)$
is identically zero almost surely. Taking this into account the last
integral in (\ref{MI}) is almost surely zero. Stochastic integration
by parts 
now gives us on $\{t<\sigma^\xi\}$
\begin{eqnarray}
 e^{-qt}v(\overline{U}^\xi_t)-v(\overline{U}^\xi_0) &=& \Lambda^\xi_t + \int_0^t  e^{-qs}(\Gamma-q)
 v(\overline{U}^\xi_{s-})ds \notag\\
 && +
\sum_{s\leq t}\mathbf{1}_{\{\Delta \overline{L}^{\xi}_s>0\}}e^{-qs}
\left\{ v(\overline{U}^\xi_{s-} + \Delta X_s - \Delta \overline{L}^\xi_s) - v(\overline{U}^\xi_{s-} +\Delta X_s) \right\}\notag\\
&&- \int_{(0,t]} e^{-qs}v^{\prime}(\overline{U}^\xi_{s-})d \overline{L}^{\xi,
c}_s, \label{IBP}
\end{eqnarray}
where $\Lambda^\xi_t = \int_0^t e^{-qs}dM^\xi_s$ is a local
martingale.

Now note that by inspection, using the properties of $a^*$, we see $v^{\prime}(x)\geq 1$ and moreover, 
on $\{\Delta \overline{L}^\xi_s>0\}$, 
\begin{equation*}
v(\overline{U}^\xi_{s-} +\Delta X_s- \Delta \overline{L}^\xi_s) -
v(\overline{U}^\xi_{s-}+\Delta X_s)= -\int_{\overline{U}^\xi_{s-} +\Delta X_s - \Delta
\overline{L}^\xi_s}^{\overline{U}^\xi_{s-}  +\Delta X_s}v^{\prime}(x)dx \leq - \Delta
\overline{L}^\xi_s.
\label{note1}
\end{equation*}
Note also that 
\begin{equation*}
\int_{(0,t]} e^{-qs} d \overline{L}^{\xi}_s = \int_{[0,t]} e^{-qs} d L^{\xi}_s -L_{0+}^\xi
\label{note2}
\end{equation*}
and by the mean value theorem and the fact that $v'(x)\geq 1$ we also have under $\mathbb{P}_x$ that $v(\overline{U}^\xi_0)\leq v(x) - L^\xi_{0+}$.
Recalling the property that $(\Gamma-q)v\leq
0$ for all $x>0$ it follows that for any appropriate localization
sequence of stopping times $\{T_n : n\geq 1\}$ we have under $\mathbb{P}_x$
\begin{eqnarray*}
v(x) - L_{0+}^\xi
&\geq&-\Lambda^\xi_{ \sigma^\xi\wedge T_n} +\int_{(0, \sigma^\xi\wedge T_n]}
e^{-qs}d\overline{L}^\xi_s+e^{-q( \sigma^\xi\wedge T_n)}
v(\overline{U}^\xi_{ \sigma^\xi\wedge T_n})\\
&\geq&-\Lambda^\xi_{ \sigma^\xi\wedge T_n} +
 \int_{[0,\sigma^\xi\wedge T_n]}
e^{-qs}dL^\xi_s - L_{0+}^\xi.
\end{eqnarray*}

 Taking expectation and then  limits as $n\uparrow\infty$ and recalling that $\xi$ is an
arbitrary strategy in $\Xi$, we thus deduce that
\[
v(x)
 \geq
 \sup_{\xi\in\Xi}
 \mathbb{E}_x\left(\int_{[0,\sigma^\xi]} e^{-qt} dL^{\xi}_t\right)
 =v^*(x).
\]
On the other hand, thanks to the
expression
$$
v(x):= \mathbb{E}_x\left(\int_{[0,\sigma^{a^*}]} e^{-qt} dL^{a^*}_t\right),
$$
the upper bound is attained by the barrier strategy at $a^*$ and the
proof is complete.
\end{proof}

\section*{Appendix: Proof of Lemmas~\ref{teoWH} and \ref{logconv}}

\begin{proof}[Proof of Lemma \ref{teoWH}]
We recall that the right continuous inverse of the local time at $0$
for $X$ reflected at its supremum, $L^{-1},$ and that of $X$
reflected  at its infimum, say $\widehat{L}^{-1},$ are possibly
killed subordinators whose Laplace exponents are given by
$\Phi(\cdot),$ and $\widehat{\kappa}(\cdot,0)$ respectively. It
follows by the time-space Wiener-Hopf factorization that
$$
q=\Phi(q)\widehat{\kappa}(q,0),\qquad q \geq 0,
$$
see e.g. \cite{Be1996} Chapter VII. Thus
$\widehat{\kappa}(q,0)=q/\Phi(q),$ $q\geq 0.$ We recall that
$\widehat{\kappa}(\cdot,\cdot)$ is the Laplace exponent of the
bivariate descending ladder subordinator, and hence it can be
represented as
\begin{equation*}
\begin{split}
\widehat{\kappa}(\lambda,\beta)=\kappa(\lambda,0)+{\rm d}_{1}\beta +
\int_{(0,\infty)^2}\mu_{-}(dt,dh)\left(e^{-\lambda t}-e^{-\lambda
t-\beta h}\right),\qquad \beta,\lambda\geq 0,
\end{split}
\end{equation*}
where ${\rm d}_{1}\geq 0$ and $\mu_{-} $ is the L\'evy measure of
the bivariate descending ladder subordinator. It follows that for $q\geq 0$ fixed $\widehat{\kappa}(q,\cdot)$ is a Bernstein
function. Since $\widehat{\kappa}(0,\beta)=\phi(\beta)$ for
$\beta\geq 0$ and the formula in the last display holds for every
$\lambda\geq 0,$ we get that ${\rm d}_{1}=\D.$ That is, the drift
term of the Bernstein function $\widehat{\kappa}(q,\cdot)$ is equal
to $\D.$ Moreover, it has been proved in Corollary 6 in \cite{DK2006}
that the measure $\mu_{-}$ can be written as
$$
\mu_{-}(dt,dh)=\int_{[0,\infty)}\mathcal{U}_{+}(dt, ds)\Pi(dh+s),
\qquad t,h>0,
$$
where $\mathcal{U}_{+}$ denotes the potential measure of the
ascending ladder subordinator. In our case $X$ is spectrally
negative and hence due to the absence of positive jumps this formula
becomes
$$
\mu_{-}(dt,dh)=\int_{[0,\infty)}\mathcal{U}_{+}(dt, ds)\Pi(dh+s) =
\int_{[0,\infty)}ds\p\left(L^{-1}_{s}\in dt\right)\Pi(dh+s),
$$
see e.g. Exercise 7.5 in \cite{Kyp}. This allows us to write
\begin{equation*}
\begin{split}
\int_{(0,\infty)^2}\mu_{-}(dt,dh)\left(e^{-q t} - e^{-qt -\beta h}
\right) &=\iiint_{(0,\infty)^3}ds\p\left(L^{-1}_{s}\in dt\right)
\Pi(dh+s)\left(e^{-q t}-e^{-qt -\beta h}\right)\\
&=\iint_{(0,\infty)\times(0,\infty)}ds\Pi(dh+s)\left( e^{-s\Phi(q)}
- e^{-s\Phi(q) - \beta h} \right)\\
&=\int^{\infty}_0(1-e^{-\beta h})\int^{\infty}_0e^{-s\Phi(q)}
\Pi(dh+s)ds, \qquad \beta\geq 0.
\end{split}
\end{equation*}
As a consequence, for $q\geq 0$ fixed, the tail of the L\'evy
measure of $\widehat\kappa(q,\cdot)$ is given by
$$
\Upsilon_{q}(z,\infty):=\int^{\infty}_{z}\int^{\infty}_0e^{-s\Phi(q)}
\Pi(dh+s)ds = e^{\Phi(q)z} \int^{\infty}_{z} du e^{-\Phi(q)u}
\Pi(u,\infty),\qquad z>0.
$$
This proves the claim about the tail of the L\'evy measure. Using
it we get that the L\'evy measure of $\widehat{\kappa}(q,\cdot)$
has a density given by
$$
\upsilon_{q}(x):=\overline{\Pi}(x)-\Phi(q)e^{\Phi(q)x}\int^\infty_{x}
e^{-\Phi(q) y}\overline{\Pi}(y)\mathrm{d}y,\qquad x> 0.
$$
To prove that $\upsilon_{q}$ is non-increasing we observe first that an integration by parts leads to the equality
\begin{equation}\label{qdens}
\begin{split}
\upsilon_{q}(x)&=\overline{\Pi}(x)-\left(\overline{\Pi}(x)-e^{\Phi(q)x}\int^\infty_{x}e^{-\Phi(q)z}\pi(z)dz\right)=e^{\Phi(q)x}\int^\infty_{x}e^{-\Phi(q)z}\pi(z)dz,
\end{split}
\end{equation}
for $x>0$.
Owing to the fact that $\pi$ is  non-increasing, thanks to the assumption that it is a log convex L\'evy density, we have that for $0<x<y$
\begin{equation*}
\begin{split}
\upsilon_{q}(x)-\upsilon_{q}(y)&=\Phi(q)e^{\Phi(q)x}\int^y_{x}e^{-\Phi(q)z}\pi(z)dz+\left(e^{\Phi(q)x}-e^{\Phi(q)y}\right)\Phi(q)\int^\infty_{y}e^{-\Phi(q)z}\pi(z)dz\\
&\geq \pi(y)e^{\Phi(q)x}\left(e^{-\Phi(q)x}-e^{-\Phi(q)y}\right)+\pi(y)\left(e^{\Phi(q)x}-e^{\Phi(q)y}\right)e^{-\Phi(q)y}=0,
\end{split}
\end{equation*} that is, $\upsilon_{q}$ is non-increasing.
\end{proof}

\begin{proof}[Proof of Lemma~\ref{logconv}]By assumption $\pi$ is  a log convex L\'evy density  (and hence non-increasing) and then by the first paragraph
of the proof of Theorem 2 in \cite{Gri78} we know that for $\beta\geq 0,$  the functions $x\mapsto e^{-\beta x}\pi(x),$ and $\int^\infty_{x}e^{-\beta z}\pi(z)dz,$ for $x>0,$ are log convex. Hence, by taking $\beta=0$ we prove the claim about $\overline{\Pi}$ in Theorem~\ref{thm:2.3}. Furthermore, it then follows that the function $x\mapsto e^{\beta x}\int^\infty_{x}e^{-\beta z}\pi(z)dz,$ for $x>0,$ is log convex. We deduce the result claimed in Lemma~\ref{logconv} by taking $\beta=\Phi(q)$ and using the characterization of the L\'evy density of $\widehat{\kappa}(q,\cdot)$ obtained in (\ref{qdens}).
\end{proof}




\section*{Acknowledgments}
The authors gratefully acknowledge funding from EPSRC
grant numbers EP/E047025/1 and  EP/C500229/1 and Royal Society grant number RE-MA1004. We are also grateful to
Ronnie Loeffen for his comments on earlier drafts of this paper, in particular with regard to the computations around (\ref{MI}).

\vspace{.5in}
\begin{singlespace}

\end{singlespace}

\end{doublespace}
\end{document}